\documentclass{gtart_h}  

\def\ifplaintex{\expandafter\ifx\csname documentclass\endcsname\relax}

\ifplaintex 
\else
\fi

\expandafter\ifx\csname beginpicture\endcsname\relax
\expandafter\ifx\csname documentclass\endcsname\relax
\input pictex \else
\input prepictex \input pictex \input postpictex \fi\fi

\def\gt{{\mathsurround=0pt\it $\cal G\mskip-2mu$eometry \&\ 
$\cal T\!\!$opology}}        

\def\gtp{{\mathsurround=0pt\it $\cal G\mskip-2mu$eometry \&\ 
$\cal T\!\!$opology $\cal P\!$ublications}}  

\def\lognumber#1{\def\thelognumber{#1}}
\def\volumenumber#1{\def\thevolumenumber{#1}}
\def\papernumber#1{\def\thepapernumber{#1}}
\def\volumeyear#1{\def\thevolumeyear{#1}}

\def\pagenumbers#1#2{\def\startpage{#1}\def\finishpage{#2}}
\def\published#1{\def\publishdate{#1}}
\def\proposed#1{\def\theproposer{#1}}
\def\seconded#1{\def\theseconders{#1}}
\def\received#1{\def\receiveddate{#1}}
\def\revised#1{\def\reviseddate{#1}}
\def\accepted#1{\def\accepteddate{#1}}

\long\def\asciiabstract#1{\long\def\theasciiabstract{#1}}

\let\\\par\let\thelognumber\relax
\let\thevolumenumber\relax\let\thepapernumber\relax
\let\thevolumeyear\relax\let\thesamplenumber\relax\let\startpage\relax
\let\finishpage\relax\let\publishdate\relax\let\receiveddate\relax
\let\reviseddate\relax\let\accepteddate\relax\let\theasciititle\relax
\let\theasciiauthors\relax
\let\theasciiabstract\relax
\let\theasciiemail\relax\let\theshortauthors\relax\let\theshorttitle\relax

\long\def\maketitlep{   

\count0=\startpage

\gt\hfill      
\beginpicture
\setcoordinatesystem units <0.33truein, 0.33truein> point at 2.2 0.9
\setplotsymbol ({$\cal G$})
\plotsymbolspacing=9truept
\circulararc 315 degrees from 0 1 center at 0 0
\setplotsymbol ({$\cal T$})
\circulararc 315 degrees from 1 -1 center at 1 0
\endpicture
\break
{\small\ifx\thesamplenumber\relax 
Volume \else Sample
\fi\thevolumenumber\ (\thevolumeyear)
\startpage--\finishpage\nl
Published: \publishdate}
\vglue 0.5truein plus 0.4fil minus 0.1truein

{\parskip=0pt\leftskip 0pt plus 1fil\def\\{\par\smallskip}{\ifplaintex\large
\else\Large\fi\bf\thetitle}\par\medskip}   

\vglue 0pt plus 0.1fil 

{\parskip=0pt\leftskip 0pt plus 1fil\def\\{\par}{\sc\theauthors}
\par\medskip}

\vglue 0pt plus 0.1fil 

{\small\parskip=0pt\let\newline\\
{\leftskip 0pt plus 1fil\def\\{\par}{\sl\theaddress}\par}
\expandafter\ifx\theemail\relax    
\relax\else\vglue 5pt plus 0.02fil minus 2pt\def\\{\stdspace{\rm 
and}\stdspace} 
\cl{Email:\stdspace\tt\theemail}\fi
\ifx\theurl\relax                  
\relax\else\vglue 5pt plus 0.02fil minus 2pt\def\\{\stdspace{\rm 
and}\stdspace}
\cl{URL:\stdspace\tt\theurl}\fi\par}

\vglue 7pt plus 0.3fil minus 3pt

{\bf Abstract}
\vglue 5pt plus 0.1fil minus 2pt

\theabstract

\vglue 7pt plus 0.3fil minus 3pt

{\bf AMS Classification numbers}\quad Primary:\quad \theprimaryclass

Secondary:\quad \thesecondaryclass

\vglue 5pt plus 0.3fil minus 2pt

{\bf Keywords:}\quad \thekeywords

\vglue 10pt plus 0.5fil minus 5pt

{\small  Proposed: \theproposer\hfill Received: \receiveddate\nl
Seconded: \theseconders\hfill 
\ifx\reviseddate\relax                         
Accepted: \accepteddate                        
\else
Revised: \reviseddate                          
\fi}
\eject
}       

\let\maketitlepage\maketitlep
\let\maketitle\maketitlepage

\font\phead=cmsl9 scaled 950
\font=cmsl9 scaled 1050
\font\pnum=cmbx10 scaled 913
\font\lnum=cmbx10 
\font\pfoot=cmsl9 scaled 950
\font=cmsl9 scaled 1050
\ifplaintex
\headline{\vbox to 0pt{\vskip -4.5mm\line{\small\phead\ifnum
\count0=\startpage ISSN 1364-0380 (on line)
1465-3060 (printed) \hfill {\pnum\folio}\else\ifodd\count0\def\\{ }%
\ifx\theshorttitle\relax\thetitle\else\theshorttitle\fi\hfill{\pnum\folio}
\else\def\\{ and }{\pnum\folio}\hfill\ifx\theshortauthors\relax\theauthors
\else\theshortauthors\fi\fi\fi}\vss}}
\footline{\vbox to 0pt{\vglue 0mm\line{\small\pfoot\ifnum\count0=\startpage
\copyright\ \gtp\hfill\else
\gt, Volume \thevolumenumber\ (\thevolumeyear)\hfill\fi}\vss
}}
\else
\makeatletter
\def\@oddhead{{\small\lhead\ifnum\count0=\startpage ISSN 1364-0380 (on line)
1465-3060 (printed) \hfill {\lnum\number\count0}\else\ifodd\count0
\def\\{ }\ifx\theshorttitle\relax \thetitle \else\theshorttitle\fi\hfill
{\lnum\number\count0}\else\def\\{ and }{\lnum\number\count0}
\hfill\ifx\theshortauthors\relax 
\theauthors\else\theshortauthors\fi\fi\fi}}\def\@evenhead{\@oddhead}
\def\@oddfoot{\small\lfoot\ifnum\count0=\startpage\copyright\ \gtp\hfill\else
\gt, Volume \thevolumenumber\ (\thevolumeyear)\hfill\fi}
\def\@evenfoot{\@oddfoot}
\makeatother
\fi

\newwrite\gtoutfile
\long\gdef\makeheadfile{  
{\def\\{, }\def\s{ }
\immediate\openout\gtoutfile head.xxx
\immediate\write\gtoutfile{Proxy-for: \ifx\theasciiauthors\relax
\theauthors\else\theasciiauthors\fi\s<\ifx\theasciiemail\relax\theemail\else\theasciiemail\fi>}
\immediate\write\gtoutfile{\noexpand\\}
\immediate\write\gtoutfile{Authors: \ifx\theasciiauthors\relax
\theauthors\else\theasciiauthors\fi}
{\def\\{ }\immediate\write\gtoutfile{Title: \ifx\theasciititle\relax
\thetitle\else\theasciititle\fi}}
\immediate\write\gtoutfile{Subj-class: GT or SG or MG etc}
\immediate\write\gtoutfile{MSC-class: \theprimaryclass\ifx\thesecondaryclass\relax\else, \thesecondaryclass\fi}
\immediate\write\gtoutfile{Journal-ref: Geom. Topol. \thevolumenumber
(\thevolumeyear) \startpage-\finishpage}
\immediate\write\gtoutfile{Comments: Published by Geometry and Topology at}
\immediate\write\gtoutfile{\s\s http://www.maths.warwick.ac.uk/gt/GTVol\thevolumenumber/paper\thepapernumber.abs.html}
\immediate\write\gtoutfile{\noexpand\\}
\immediate\write\gtoutfile{}
\ifx\theasciiabstract\relax
\immediate\write\gtoutfile{\theabstract}\else
\immediate\write\gtoutfile{\theasciiabstract}\fi
\immediate\write\gtoutfile{}
\immediate\write\gtoutfile{\noexpand\\}
\immediate\write\gtoutfile{}
\immediate\closeout\gtoutfile}}  

\def\maketitlepage{\maketitlep\makeheadfile}
\let\maketitle\maketitlepage

\lognumber{385}
\volumenumber{8}\papernumber{16}\volumeyear{2004}
\pagenumbers{645}{673}
\received{25 November 2003}
\revised{21 April 2004}
\published{22 April 2004}

\accepted{13 March 2004}

\proposed{Thomas Goodwillie}
\seconded{Ralph Cohen, Haynes Miller}

\usepackage{amssymb, amsmath,amscd}
\usepackage{verbatim}

\newtheorem{theorem}{Theorem}[section]
\newtheorem{lemma}[theorem]{Lemma}
\newtheorem{proposition}[theorem]{Proposition}
\newtheorem{corollary}[theorem]{Corollary}

\theoremstyle{definition}
\newtheorem{definition}[theorem]{Definition}

 \theoremstyle{remark}
\newtheorem{remark}[theorem]{Remark}
\newtheorem{example}[theorem]{Example}

\makeatletter\let\c@equation=\c@theorem\makeatother

\numberwithin{equation}{section}

\newcommand{\cy}{\text{cy}}
\newcommand{\colim}{\operatornamewithlimits{colim}}
\newcommand{\D}{\mathcal D}

\newcommand{\HH}{\operatorname{HH}}
\newcommand{\I}{\mathcal I}

\newcommand{\hocolim}{\operatornamewithlimits{hocolim}}
\newcommand{\K}{\operatorname{K}}
\newcommand{\lra}{\longrightarrow}
\newcommand{\lla}{\longleftarrow}
\newcommand{\la}{\leftarrow}

\newcommand{\Map}{\operatorname{Map}}

\newcommand{\Th}{\operatorname{TH}}
\newcommand{\THH}{\operatorname{THH}}
\newcommand{\tr}{\text{tr}}

\begin{document}

\title{Units of ring spectra and their traces\\in algebraic K-theory}
\author{Christian Schlichtkrull}
\address{Department of Mathematics, Oslo University\\PO Box 1053, 
Blindern\\NO-0316 Oslo, Norway}
\email{krull@math.uio.no}

\begin{abstract}
Let $GL_1(R)$ be the units of a commutative ring spectrum $R$. In this
paper we identify the composition
$$
\eta_R\co BGL_1(R)\to
\K(R)\stackrel{}{\to}\THH(R)\stackrel{}{\to}
\Omega^{\infty}(R),
$$
where $\K(R)$ is the algebraic K-theory and $\THH(R)$ the
topological Hochschild homology of $R$. As a corollary we show
that classes in $\pi_{i-1}R$ not annihilated by the stable Hopf map
$\eta\in \pi_1^s(S^0)$ give rise to non-trivial classes in
$\K_i(R)$ for $i\geq 3$.
\end{abstract}

\asciiabstract{%
Let GL_1(R) be the units of a commutative ring spectrum R. In this
paper we identify the composition

BGL_1(R)->K(R)->THH(R)->\Omega^{\infty}(R),

where K(R) is the algebraic K-theory and THH(R) the topological
Hochschild homology of R. As a corollary we show that classes in
\pi_{i-1}(R) not annihilated by the stable Hopf map give rise to
non-trivial classes in K_i(R) for i\geq 3.}




\maketitlepage

\section{Introduction}
Given a connective (symmetric) ring spectrum $R$, we follow Waldhausen
and define the  
units $GL_1(R)$ to be the union of the components in $\Omega^{\infty}(R)$
that correspond to units in the discrete ring
$\pi_0R$. With this definition $GL_1(R)$ is a
group-like monoid whose group of components equals $GL_1(\pi_0R)$. As
in the case of a discrete ring there is a natural map $BGL_1(R)\to
\K(R)$ to the algebraic K-theory of $R$. If $R$ is a commutative
discrete ring this is split by the determinant, but the
definition of the determinant does not generalize to the setting of
ring spectra and the above map is in general not split,
even if $R$ is commutative. For example, Waldhausen shows \cite{W2}
that this fails quite badly for the sphere spectrum. However, it
turns out that the notion of traces of matrices does generalize to
ring spectra. This gives rise to the trace map $\tr\co \K(R)\to \THH(R)$,
where the target is the topological Hochschild homology first defined
by B\"okstedt \cite{Boe}.
The purpose of the present paper is to identify the composition
\begin{equation}\label{compositionequation}
\eta_R\co BGL_1(R)\to
\K(R)\stackrel{\tr}{\to}\THH(R)\stackrel{r}{\to
}
\Omega^{\infty}(R)
\end{equation}
when $R$ is a commutative ring spectrum. The first two arrows are
defined for any (symmetric) ring spectrum, whereas the definition
of the last map depends on $R$ being commutative. By definition,
$\THH(R)$ is the infinite loop space associated to the realization of
the cyclic spectrum $[k]\mapsto R^{\wedge(k+1)}$ with Hochschild type
structure maps. We shall use B\"okstedt's explicit definition of the
smash products $R^{\wedge(k+1)}$. If $R$ is
commutative, the degree-wise
multiplication $R^{\wedge(k+1)}\to R$ defines a map to the constant
cyclic spectrum. This gives rise to the infinite loop map $r$ in the
definition of $\eta_R$. 

In order to state our main result, we need the fact that $GL_1(R)$
has the structure of an infinite loop space when $R$ is
commutative, ie, that there exists a spectrum $gl_1(R)$ such that
$\Omega^{\infty}(gl_1(R))\simeq GL_1(R)$. (We follow the
convention to use small letters for the
spectrum associated to an infinite loop space written in capital
letters.) It will be convenient for our purpose to give an explicit
construction of $gl_1(R)$ using Segal's notion of $\Gamma$-spaces.
Let $\eta\in \pi^s_1(S^0)$ denote the stable Hopf map.
\begin{theorem}\label{introductiontheorem}
The composite map $\eta_R$ admits a factorization
$$
BGL_1(R)\to GL_1(R)\to \Omega^{\infty}(R),
$$
in which the second map is the natural inclusion and the first map is
multiplication by $\eta$ in the sense of the following 
commutative diagram in the homotopy category of spaces,
$$
\begin{CD}
BGL_1(R) @>>> GL_1(R) \\
@VV\sim V @VV \sim V \\
\Omega^{\infty}(gl_1(R)\wedge S^1) @>\Omega^{\infty}(\text{id}\wedge\eta)
>>\Omega^{\infty}(gl_1(R)).
\end{CD}
$$
\end{theorem}
In the case where $R$ equals the sphere spectrum
this result is due to  B\"okstedt and Waldhausen \cite{BW} (with a
completely different proof).

It is clear from the definition that there is an isomorphism of
abelian groups $\pi_igl_1(R)\cong \pi_iR$ for
$i\geq 1$, but since the spectrum structures are different this is not an
isomorphism of $\pi_*^s(S^0)$-modules. However, using that $\eta$ is
realized as an unstable map $\eta\co S^3\to S^2$, it is not difficult to
check that the actions of $\eta$ are compatible in degrees $i\geq
2$. The following is then an immediate corollary of Theorem
\ref{introductiontheorem}.
\begin{corollary}
For $i\geq 3$, the composition
$$
\pi_{i-1}R\cong\pi_iBGL_1(R)\to\pi_i\K(R)\to\pi_i\THH(R)\to\pi_i R
$$
is multiplication by $\eta\in \pi_1^s(S^0)$.
\end{corollary}
It thus follows that classes in $\pi_{i-1}R$ not annihilated
by $\eta$ give rise to non-trivial elements in $\pi_i\K(R)$.
\begin{example}
Let $R=ko$, the real connective K-theory spectrum. In this case
$GL_1(ko)\simeq\{\pm1\}\times BO_{\otimes}$, where $\otimes$
  indicates that the
H-space structure is the one corresponding to tensor products of
vector bundles. Using the cofibration sequence $\Sigma
ko\stackrel{\eta}{\to}ko\to ku$, \cite[V.5.15]{MQR}, we see that
$$
\mathbb Z\cong\pi_{8k}(ko)\stackrel{\eta}{\to}
\pi_{8k+1}(ko)\cong\mathbb Z/2
$$
is surjective and that
$$
\mathbb
Z/2\cong\pi_{8k+1}(ko)\stackrel{\eta}{\to}\pi_{8k+2}(ko)\cong\mathbb
Z/2
$$
is an isomorphism. We conclude that for $k\geq 1$,
\begin{itemize}
\item
$\pi_{8k+1}BBO_{\otimes}\cong\mathbb Z$ maps non-trivially to
$\pi_{8k+1}\K(ko)$;

\medskip
\item
$\pi_{8k+2}BBO_{\otimes}\cong \mathbb Z/2$ injects as a direct
summand in $\pi_{8k+2}\K(ko)$.
\end{itemize}
\end{example}

This example is interesting in view of the attempts \cite{AR}, 
\cite{BDR}, to relate
algebraic K-theory to elliptic cohomology and the chromatic filtration
of homotopy theory. Another major source for the interest in 
algebraic K-theory in the non-linear setting is the relation to high
dimensional manifold theory via Waldhausen's work on stable
concordances \cite{W1}.
\begin{example}
Let $R=\Sigma^{\infty}(G_+)$ be the suspension spectrum of a
commutative (or $E_{\infty}$) group-like monoid  $G$. By
definition, the algebraic K-theory of this spectrum is Waldhausen's
$A(BG)$. In this case, $\pi_iBGL_1(R)\cong \pi_{i-1}^s(G_+)$, and
thus classes in the stable homotopy that are not annihilated by $\eta$
map non-trivially to $\pi_iA(BG)$ in degrees $i\geq 3$.
\end{example}

\begin{remark}
Given a discrete ring $\bar R$, the algebraic K-theory of the
associated Eilenberg-MacLane spectrum $H\bar R$ reduces to Quillen's
$\K(\bar R)$. Starting with a ring spectrum $R$ and $\bar R=\pi_0 R$,
the linearization map $R\to H\bar R$ gives rise to a fibration
sequence
$$
F\to \K(R)\to \K(\bar R),
$$
where by definition $F$ is the homotopy fibre. Let $SL_1(R)$ be
the unit component of $GL_1(R)$. Using that $BSL_1(\bar R)=*$ we get a
map $BSL_1(R)\to F$ which is important in the
understanding of how algebraic K-theory behaves under linearization.
\end{remark}

The proof of Theorem \ref{introductiontheorem} breaks up into two
parts. The first part is to give a description of $\eta_R$ in
non-K-theoretical terms as the composition 
$$
BGL_1(R)\to L(BGL_1(R))\stackrel{\sim}{\leftarrow}
B^{cy}GL_1(R)\stackrel{r}{\to}GL_1(R)\subseteq
\Omega^{\infty}(R).
$$
Here $L(BGL_1(R))$ denotes the free loop space of $BGL_1(R)$ and
$B^{\cy}GL_1(R)$ is  Waldhausen's cyclic bar construction, see Section
\ref{K-theorysection}. The first map is the inclusion of the constant
loops and the map $r\co B^{\cy}GL_1(R)\to GL_1(R)$ is given by
iterated multiplication in $GL_1(R)$. The fact that $GL_1(R)$ is
an infinite loop space ensures that it is sufficiently homotopy
commutative for the latter map to be well-defined. 

The second part of the proof is then to show that the composite map 
$BGL_1(R)$ $\to GL_1(R)$ is multiplication by $\eta$. This
follows from a general analysis of how the free loop space of an infinite
loop space relates to the cyclic bar construction. Let us say that a
sequence of maps of based spaces $F\to X\to Y$ is a homotopy fibration
sequence if (i) the composition is constant and (ii) the canonical map
from $F$ to the homotopy fiber of the second map is a weak homotopy
equivalence. (This definition is most useful if $Y$ is connected.) 
Given a well-pointed group-like topological monoid $G$, there is a
commutative diagram of homotopy fibration sequences
$$
\begin{CD}
G@>>> B^{\cy}G @>>>BG\\
@VV\sim V @VV\sim V @|  \\
\Omega(BG)@>>> L(BG) @>>> BG,
\end{CD}
$$
in which the lower sequence is split by the inclusion of the constant
loops $BG\to L(BG)$. If furthermore $G$ admits the structure of an
infinite loop space, then the upper sequence has a natural splitting
$B^{\cy}G\to G$ given by the iterated product in $G$. The failure of
these splittings to be compatible is measured by the fact that the
composition
$$
BG\to L(BG)\simeq B^{\cy}G\to G
$$
is multiplication by $\eta$ in the sense described above for $GL_1(R)$.    

\medskip

The paper is as a whole fairly self-contained, and in particular we
present in    
Section \ref{tracesection} a new explicit construction of the
trace map $\tr\co \K(R)\to \THH(R)$. This version of the trace map is
used here to 
identify the action on $BGL_1(R)$, but there are many other
applications of this combinatorial construction. In
Section \ref{unitssection} we recall the definition of symmetric
ring spectra and their units and in Section
 \ref{K-theorysection} we recall Waldhausen's  definition of
algebraic K-theory in this framework. The Sections
\ref{unitssection}-\ref{tracesection} can be read as a self-contained
account of the topological trace map.

In Section \ref{gammasection} we explain the infinite loop structure
of $GL_1(R)$ used in the formulation of Theorem
\ref{introductiontheorem}, and in Section \ref{splittingsection} we
construct the splitting
$r\colon \THH(R)\to\Omega^{\infty}(R)$ and complete the first part of the proof.
Finally, in Section \ref{Hopfsection} we consider the relationship
between the free loop space and the cyclic bar construction of an
infinite loop space and finish the second part of the proof.
\medskip

\subsection{Notation and conventions} Let $\mathcal T$ be the
category of based spaces. In this paper this can be understood as
either the category of compactly generated Hausdorff (or weak
Hausdorff) topological spaces or the category of based simplicial
sets. However, we will usually use the topological terminology and
talk about topological monoids etc. 
In both cases equivalences mean weak homotopy equivalences.
In the topological case we will sometimes have to assume that base
points are non-degenerate in the usual sense of being neighborhood
deformation retracts.   

We let $S^n$ denote the $n$-fold smash product of the circle
$S^1=I/\partial I$. By a spectrum $E$ we understand a sequence
$\{E_n\co n\geq 0\}$ of based spaces together with based maps
$\sigma\co S^1\wedge E_n\to E_{n+1}$. Again this may be interpreted
either in the topological or simplicial category. A map of spectra
$f\co E\to F$ is a sequence of based maps $f_n\co E_n\to F_n$ that
commute with the structure maps. We say that $f$ is an equivalence
if it induces an isomorphism on spectrum homotopy groups, the
latter being defined by $\pi_nE=\colim_k\pi_{n+k}E_{k}$. All spectra
we consider will be connective, ie, $\pi_nE=0$ for $n<0$. We
shall also assume that the spectra we consider are
\emph{convergent} in the sense that there exists an unbounded, non-decreasing
sequence of natural numbers $\{\lambda_n\co n\geq 0\}$ such that
$S^1\wedge E_n\to E_{n+1}$ is as least $n+\lambda_n$-connected for
all $n$. This is not a serious restriction as any connective spectrum is
equivalent to a convergent one.

\section{Units of ring spectra}\label{unitssection}
In this section we recall Waldhausen's definition of the space of
units associated to a ring spectrum. We shall work in the
framework of symmetric spectra and begin by recalling the relevant
definitions from \cite{HSS} and, for the version with topological
spaces instead of simplicial sets, \cite{MMSS}. 

\subsection{Symmetric spectra}
A symmetric spectrum is a spectrum in
which each of the spaces $E_n$ is equipped with a base point
preserving left $\Sigma_n$-action, such that the iterated
structure maps
$$
\sigma^m\co S^m\wedge E_n\to E_{m+n}
$$
are $\Sigma_m\times\Sigma_n$-equivariant. A symmetric ring
spectrum is a symmetric spectrum equipped with
$\Sigma_n$-equivariant maps $1_n\co S^n\to E_n$ for $n\geq 0$, and
$\Sigma_m\times \Sigma_n$-equivariant maps $\mu_{m,n}\co E_m\wedge
E_n\to E_{m+n}$ for $m, n\geq 0$. In order to formulate the
axioms, let $\check\sigma^n$ be the composite
$$
\check\sigma^n\co E_m\wedge S^n\stackrel{\text{tw}}{\lra}S^n\wedge E_m
\stackrel{\sigma^n}{\lra}E_{n+m}\stackrel{\tau_{n,m}}{\lra}E_{m+n},
$$
where $\text{tw}$ twists the two factors, and
$\tau_{n,m}$ is the $(n,m)$-shuffle $i\mapsto i+m$ for $i\leq n$,
$i\mapsto i-n$ for $i>n$. Notice that $\check\sigma^n$ is
$\Sigma_m\times \Sigma_n$-equivariant.
Also, let $\sigma^0\co S^0\wedge E_n\to E_n$
and $\check\sigma^0\co E_n\wedge S^0\to E_n$ be the
canonical identifications.
These maps are required to satisfy the following relations for all
$l,m,n\geq 0$:
\begin{description}
\smallskip
\item[(a)]
$1_{m+n}=\sigma^m\circ (S^m\wedge 1_n)$,
\smallskip
\item[(b)]
$\sigma^m=\mu_{m,n}\circ (1_m\wedge E_n),\quad
\check\sigma ^n=\mu_{m,n}\circ (E_m\wedge 1_n)$,
\smallskip
\item[(c)]
$\mu_{l+m,n}\circ (\mu_{l,m}\wedge E_n)=\mu_{l,m+n}\circ (E_l\wedge\mu_{m,n})$.
\end{description}
Here condition (a) states that the maps $1_n$ assemble to give a
map of spectra $1\co S\to E$, where $S$ denotes the sphere spectrum.
Notice that (b) and (c) imply that
$$
\mu_{l,m+n}\circ(E_l\wedge\sigma^m)=\mu_{l+m,n}\circ(\check\sigma^m\wedge E_n)
$$
as maps $E_l\wedge S^m\wedge E_n\to E_{l+m+n}$ and that
$$
\sigma^l\circ (S^l\wedge \mu_{m,n})=\mu_{l+m,n}\circ(\sigma^l\wedge
E_n).
$$
These are exactly the conditions for the maps $\mu_{m,n}$ to
produce a map of spectra $\mu\co E\wedge E\to E$, where the domain is
the internal smash product in the category of symmetric spectra.
Condition (b) then says that $1$ is a two-sided unit, and (c) is
the condition that the multiplication is associative. (These comments
on the internal smash product are only to motivate the definitions; we
shall not make explicit use of the internal smash product in this
paper.) We say that $R$ is commutative if the diagrams
$$
\begin{CD}
R_m\wedge R_n@>\mu_{m,n}>>R_{m+n}\\
@VV\text{tw} V @VV\tau_{m,n} V \\
R_n\wedge R_m  @>\mu_{m,n}>>R_{n+m}
\end{CD}
$$
are commutative. 

\subsection{$\I$-spaces and $\I$-monoids}
In order to define the units of a sym\-metric ring spectrum
we need a combinatorial framework to keep track of the suspension
coordinates. Let $\I$ be the category whose objects are
the finite sets $\mathbf n=\{1,\dots,n\}$ and whose morphisms are
the injective (not necessarily order preserving) maps. The empty
set $\mathbf 0$ is an initial object. The concatenation $\mathbf
m\sqcup \mathbf n$ defined by letting $\mathbf m$ correspond to
the first $m$ elements and $\mathbf n$ to the last $n$ elements of
$\{1,\dots,m+n\}$ gives $\I$ the structure of a symmetric monoidal
category. The symmetric structure is given by the shuffles 
$\tau_{m,n}\co\mathbf m\sqcup \mathbf n\to \mathbf n\sqcup
\mathbf m$. 

We define an \emph{$\I$-space} to be a functor
$X\co\I\mathcal\to \mathcal T$. Given an 
$\I$-space $X$, we write $X_{h\I}=\hocolim_{\I}X$.
 The homotopy type of $X_{h\I}$ can
be analyzed using the following lemma due to B\"okstedt. For 
published versions see \cite[2.3.7]{M} and \cite[2.5.1]{B}. 
Let $\mathcal F_n\I$ be
the full subcategory of $\I$ containing the objects of cardinality
at least $n$. 
\begin{lemma}[B\"okstedt]\label{approximationlemma}
Let $X$ be an $\I$-space and suppose that each morphism $\mathbf
n_1\to \mathbf n_2$ in $\mathcal F_n\I$ induces a
$\lambda_n$-connected map $X(n_1)\to X(n_2)$. Then, given 
any object $\mathbf m$ in $\mathcal F_n\I$, the
natural map $X(m)\to X_{h\I}$ given by the inclusion
in the $0$-skeleton is at least $(\lambda_n-1)$-connected.
\qed
\end{lemma}
 Let us say that an
$\I$-space $X$ is \emph{convergent} if there exists an unbounded,
non-decreasing sequence of natural numbers $\{\lambda_n\co n\geq 0\}$
such that any morphism $\mathbf n_1\to \mathbf n_2$ in $\mathcal
F_n\I$ induces a $\lambda_n$-connected map $X(n_1)\to X(n_2)$. It
follows from B\"okstedt's lemma that in this case
$X_{h\I}$ is equivalent to the usual telescope of the sequence of
spaces $X(n)$ obtained by restricting to the natural subset inclusions in
$\I$. In particular, $\pi_*X_{h\I}$ is the usual directed
colimit of the groups $\pi_*X(n)$ if $X$ is convergent.
 
We say that an $\I$-space $X$ is an \emph{$\I$-monoid} if it comes
equipped with an associative and unital natural transformation
$$
\mu_{m,n}\co X(m)\times X(n)\to X(m+n),
$$
where both sides are considered functors on $\I^2$. The
unital condition means that the basepoint in $X(0)$ acts as a unit and
associativity means that the identity
$$
\mu_{l,m+n}\circ(X_l\times \mu_{m,n})=\mu_{l+m,n}\circ
(\mu_{l,m}\times X_n)
$$
holds for all $l,m,n\geq 0$. By definition an $\I$-monoid $X$ is
commutative if the diagrams
$$
\begin{CD}
X(m)\times X(n)@>\mu_{m,n} >> X(m+n)\\
@VV \text{tw} V @VV\tau_{m,n}V \\
X(n) \times X(m) @> \mu_{n,m} >> X(n+m)
\end{CD}
$$
are commutative. If $X$ is an $\I$-monoid, then $X_{h\I}$ inherits the
structure of a topological monoid. The product is given by the
composition
$$
X_{h\I}\times X_{h\I}=\hocolim_{\I\times \I}X(m)\times
X(n)\stackrel{\mu}{\to} \hocolim_{\I\times \I}X(m+n)\to X_{h\I},
$$
in which the last map is induced by the monoidal structure of $\I$. 
We say that $X$ is group-like if this is the case for
$X_{h\I}$, ie, if the monoid of components $\pi_0X_{h\I}$ is a
group. We will show in Section \ref{gammasection} that if $X$ is
commutative and group-like, then $X_{h\I}$ has the structure
of an infinite loop space.

\begin{remark}
For $\I$-spaces $X$ that are not convergent, the homotopy type of
$X_{h\I}$ may well differ from that of the usual telescope. 
Consider for example the $\I$-monoid 
$\mathbf n\mapsto B\Sigma_n$. In this case the associated homotopy
colimit is equivalent to the base point component of $Q(S^0)$. 
To see this one uses that the natural map
$B\Sigma_{\infty}\to \hocolim_{\I}B\Sigma_n$
induces an isomorphism on integral homology. By the universal
property of Quillen's plus-construction and the fact that the target
is a connected H-space, it follows that the latter is equivalent to 
$B\Sigma_{\infty}^+$. The conclusion then follows from the
Barratt-Priddy-Quillen-Segal Theorem.
As a second example, let $R$ be a discrete ring
and consider the $\I$-monoid defined by the classifying spaces
$BGL_n(R)$. By an argument similar to the above, the associated 
homotopy colimit is equivalent to the base point component of 
the algebraic K-theory
space $\K(R)$. In these examples (and many more), evaluating the
homotopy colimit over $\I$ thus has the same effect as Quillen's
plus-construction.
\end{remark}

\subsection{Units of ring spectra}\label{unitsection}
Given a symmetric ring spectrum $R$, the sequence of spaces 
$\Omega^n(R_n)$ defines an $\I$-space as follows. A  morphism
$\alpha\co\mathbf m\to \mathbf n$ in $\I$ induces a map
$\Omega^m(R_m)\to \Omega^n(R_n)$ by taking $f\in \Omega^m(R_m)$ to the
composition
\begin{equation}\label{Jfunctoriality}
S^n\stackrel{\bar\alpha^{-1}}{\lra}S^n=S^l\wedge
S^m\stackrel {S^l\wedge f}{\lra} S^l\wedge R_m\stackrel{\sigma^l}{\lra}
R_n\stackrel {\bar\alpha}{\lra}R_n.
\end{equation}
Here $\bar\alpha\co\mathbf n=\mathbf l\sqcup\mathbf m\to \mathbf n$
is the unique permutation that is order preserving on the first
$l=n-m$ elements and acts as $\alpha$ on the last $m$ elements.
The action on $S^n$ is the usual left action. The multiplication
in $R$ gives a multiplicative structure
\begin{gather*}
\mu_{m,n}\co\Omega^m(R_m)\times \Omega^n(R_n)\to
\Omega^{m+n}(R_{m+n}),\\
\mu_{m,n}(f,g)\co S^m\wedge S^n\stackrel{f\wedge
  g}{\lra} R_m\wedge R_n\stackrel{\mu_{m,n}}{\lra} R_{m+n},
\end{gather*}
which is commutative if $R$ is. We let $\Omega^n(R_n)^*$ be the
union of the components in $\Omega^n(R_n)$ that have stable
multiplicative homotopy inverses in the following sense: For each $f$ in
$\Omega^n(R_n)^*$ there exists an element $g\in \Omega^m(R_m)$ such that
$\mu_{n,m}(f,g)$ and $\mu_{m,n}(g,f)$ are homotopic to the unit
$1_{m+n}$ in $\Omega^{m+n}(R_{m+n})$. We consider
$\Omega^n(R_n)^*$ as a based space with base point $1_n$ and
restricting the above structure maps gives an
$\I$-monoid $\Omega^n(R_n)^*$. We define
$$
GL_1(R)=\hocolim_{\I}\Omega^n(R_n)^*
$$
with the monoid structure explained above. 
If $R$ is convergent so is the $\I$-space $\Omega^n(R_n)^*$, and
by Lemma \ref{approximationlemma},
$\pi_0(GL_1(R))=GL_1(\pi_0(R))$. If furthermore $R$ is commutative, the
general construction in section \ref{gammasection} will produce a
spectrum $gl_1(R)$ such that $\Omega^{\infty}(gl_1(R))\simeq GL_1(R)$.  

\section{K-theory and cyclic K-theory of ring spectra}\label{K-theorysection}
In this section we recall the definition of the algebraic K-theory
$\K(R)$ and the cyclic algebraic K-theory $\K^{\cy}(R)$ of a
symmetric ring spectrum $R$. We also recall the inclusion of the units
$BGL_1(R)\to \K(R)$. This material is due mainly to Waldhausen.
Let $M_n(R)$ be the symmetric
ring spectrum whose $m$th space is $\Map_*(\mathbf n_+,\mathbf
n_+\wedge R_m)$. The multiplication
resembles multiplication of $n\times n$ matrices over an ordinary
ring. (In this case the ``matrices'' in question have at most one
non-base point entry in each column.) We let
$GL_n(R)=GL_1(M_n(R))$ with the monoid structure coming from the
multiplication in $M_n(R)$. Using the natural maps
\begin{align*}
&\Map_*(\mathbf m_+\wedge S^k,\mathbf m_+\wedge R_k)\times
\Map_*(\mathbf n_+\wedge S^l,\mathbf n_+\wedge R_l)\\
&\to \Map_*((\mathbf m\sqcup \mathbf n)_+\wedge S^{k+l},(\mathbf
m\sqcup \mathbf n)_+\wedge R_{k+l})
\end{align*}
we have a notion of block sum of matrices and corresponding
monoid homomorphisms
$$
GL_m(R)\times GL_n(R)\to GL_{m+n}(R).
$$
These homomorphisms are associative in the obvious sense and thus
the induced maps of classifying spaces give $\coprod_{n\geq
0}BGL_n(R)$ the structure of an associative topological monoid. By
definition $\K(R)$ is the group completion
$$
\K(R)=\Omega B\bigg(\coprod_{n\geq 0}BGL_n(R)\bigg).
$$
Notice that this is the version of
algebraic K-theory with $\pi_0\K(R)=\mathbb Z$. 

The classifying space
of the units $BGL_1(R)$ embeds in the $1$-simplices of
$B_{\bullet}(\coprod_{n\geq 0}BGL_n(R))$ and since there is just
a single $0$-simplex there is an induced map
$$
S^1\wedge BGL_1(R)_+\to B\bigg(\coprod_{n\geq 0}BGL_n(R)\bigg)
$$
whose adjoint is the requested map $BGL_1(R)\to \K(R)$. The image
is contained in the $1$-component of $\K(R)$.

There is a variant of all this using the cyclic bar construction
$B^{\cy}GL_n(R)$. Recall that for a topological monoid $G$,
$B^{\cy}G$ is the realization of the cyclic space $[k]\mapsto
G^{k+1}$ with simplicial operators
\begin{align*}
&d_i(g_0,\dots,g_k)=
\begin{cases}(g_0,\dots,g_ig_{i+1},\dots,g_k),&\text{for }0\leq i<k \\
(g_{k}g_0,\dots,g_{k-1}),& \text{for } i=k
\end{cases}\\
&s_i(g_0,\dots, g_k)=(g_0,\dots,g_i,1,\dots,g_k),\quad\text{for
}0\leq i \leq k,
\end{align*}
and cyclic operator
$t_k(g_0,\dots,g_k)=(g_k,g_0,\dots,g_{k-1})$. We refer the reader to
\cite{J} for background material on cyclic spaces.
The degree-wise projections $(g_0,\dots,g_k)\mapsto (g_1,\dots,g_k)$
define a simplicial map
$p\co B^{\cy}_{\bullet}G\to B_{\bullet}G$, and if $G$ is group-like and has
a non-degenerate unit there results a homotopy fibration sequence
$$
G\to B^{\cy}G\to BG.
$$
Since $B^{\cy}_{\bullet}G$ is a cyclic space its realization has a
canonical action of the circle group $\mathbb T$.
Consider the composite map $\mathbb T\times B^{\cy}G\to
B^{\cy}G\stackrel{p}{\to} BG$, where the first map is given by the
$\mathbb T$-action. Letting $L(-)$ denote the free loop space, the
adjoint is a map $B^{\cy}G\to  L(BG)$. It is immediate from the
definition that this is a $\mathbb T$-equivariant map when the
action on $L(BG)$ is by multiplication in $\mathbb T$. The
following proposition is well-known and follows easily from the
definition of the $\mathbb T$-action. We shall prove a related
result in Proposition \ref{gammafreeloopproposition}
with a proof that can easily be adapted to the present situation.
\begin{proposition}\label{cyclicloopproposition}
There is a commutative diagram
$$\label{cyclicdiagram}
\begin{CD}
G @>i>> B^{\cy}G @>p>> BG \\
@VVV @VVV @| \\
\Omega(BG) @>i>> L(BG)@>\text{ev} >> BG
\end{CD}
$$
in which the lower sequence is the usual fibration sequence
associated to the evaluation at the unit element of $\mathbb T$.
If $G$ is group-like and has a non-degenerate unit,
then the upper sequence is a homotopy fibration sequence and the
vertical maps are equivalences.\qed
\end{proposition}
Notice that the lower sequence in (\ref{cyclicdiagram}) is split
by the inclusion of $BG$ in $L(BG)$ as the constant loops.
By definition the cyclic K-theory of a symmetric ring spectrum $R$ 
is given by
$$
\K^{\cy}(R)=\Omega B\bigg(\coprod_{n\geq 0}B^{\cy}GL_n(R)\bigg).
$$
The projections $p\co B^{\cy}GL_n(R)\to BGL(R)$ induce a map
$p\co\K^{\cy}(R)\to \K(R)$ which has a section in the homotopy
category. The quickest way to see this is to consider the diagram
of monoid homomorphisms
$$ \coprod_{n\geq
0}BGL_n(R)\lra\coprod_{n\geq
0}L(BGL_n(R))\stackrel{\sim}{\lla}\coprod_{n\geq 0}B^{\cy}GL_n(R),
$$
where the equivalence is a consequence of Proposition
\ref{cyclicloopproposition}. After group completion there results
a well-defined homotopy class $\K(R)\to\K^{cy}(R)$, giving a
section of $p$ up to homotopy.

\section{Topological Hochschild homology and the trace
map}\label{tracesection}
In this section we present an
explicit construction of the topological trace map $\tr\co\K(R)\to \THH(R)$,
where the target is the topological Hochschild homology. In order
to motivate the construction we first recall the linear trace map
for an ordinary discrete ring $\tr\co \K(R)\to \HH(R)$ with target the
Hochschild homology of $R$. The latter is the realization of a
cyclic abelian group $\HH_{\bullet}(R)\co[k]\mapsto R^{\otimes k+1}$
with cyclic structure maps similar to $B^{\cy}_{\bullet}G$. The
multi-trace $\tr\co\HH_{\bullet}(M_n(R))\to \HH_{\bullet}(R)$ is the
cyclic map given in degree $k$ by
$$
\tr(A^0\otimes\dots\otimes
A^k)=\sum_{s_0,\dots,s_k}a^0_{s_k,s_0}\otimes \dots \otimes
a^k_{s_{k-1},s_k},\quad \text{where }A^i=(a^i_{s,t}).
$$
Composing with the obvious inclusions $B^{\cy}_{\bullet}GL_n(R)\to
\HH_{\bullet}(M_n(R))$ we get a cyclic map
$$
\coprod_{n\geq 0}B_{\bullet}^{\cy}GL_n(R) \to \coprod_{n\geq
0}\HH_{\bullet}(M_n(R)) \to \HH_{\bullet}(R).
$$
This is a monoid homomorphism with respect to block-sums of
matrices on the domain and the abelian group structure on the
target. After realization and group completion we get maps
$$
\K^{\cy}(R)\to \Omega B\HH(R)\stackrel{\sim}{\la}\HH(R).
$$
The linear trace map $\tr\co\K(R)\to\HH(R)$ is the homotopy class
obtained by composing with the homotopy section $\K(R)\to
\K^{\cy}(R)$.

\subsection{Topological Hochschild homology}
Topological Hochschild homology is obtained by replacing the
tensor products in $\HH_{\bullet}(R)$ by smash products of
spectra. We shall follow B\"okstedt \cite{Boe} 
in making this precise. Given
a sequence of symmetric spectra $E^1,\dots, E^r$, we consider
the smash product as a multi-indexed spectrum in the natural
way, 
$$
(E^1\wedge\dots\wedge E^r)_{n_1,\dots,n_r}=E^1_{n_1}\wedge\dots
\wedge E^r_{n_r}.
$$
In general an $r$-fold multi-indexed symmetric spectrum
$E=\{E_{n_1,\dots,n_r}\}$ has an associated infinite
loop space
$$
\Omega^{\infty}(E)=\hocolim_{\I^r}\Omega^{n_1+\dots+n_r}(E_{n_1,\dots,n_r}).
$$
The functoriality underlying this definition is analogous to
that in (\ref{Jfunctoriality}). We shall always use the symbol
$\Omega^{\infty}(E)$ in this precise way. Notice that the monoidal
structure of $\I^r$ makes $\Omega^{\infty}$ a functor from
multi-indexed spectra to topological monoids. The topological
Hochschild homology of a ring spectrum $R$ is the topological 
realization of
the cyclic spectrum $\Th_{\bullet}(R)$, defined in spectrum degree
$n$ by
$$
\Th_{\bullet}(R,n)\co[k]\mapsto \Omega^{\infty}(\underbrace{R\wedge\dots\wedge
R}_{k+1}\wedge S^n).
$$
The spectrum structure maps are defined in the obvious way
involving only the $S^n$-factor. This construction represents
B\"okstedt's solution to the problem of how to turn the
multi-indexed spectrum $R^{\wedge(k+1)}$ into an equivalent
singly-indexed spectrum. The cyclic structure maps are analogous
to those in $B_{\bullet}^{\cy}(G)$ and $\HH_{\bullet}(R)$. Thus
for example $d_0\co\Th_1(R)\to \Th_0(R)$ is the composition
\begin{align*}
\hocolim_{\I^2}\Omega^{n_0+n_1}(R_{n_0}\wedge R_{n_1}\wedge
S^n)&\to\hocolim_{\I^2}\Omega^{n_0+n_1}(R_{n_0+n_1}\wedge
S^n)\\
&\to\hocolim_{\I}\Omega^{n_0}(R_{n_0}),
\end{align*}
where the first map uses the multiplication in $R$ and the second
map is induced by the monoidal structure $\sqcup\co\I\times \I\to \I$.
It follows from the version of B\"okstedt's approximation Lemma
\ref{approximationlemma} with $\I^{k+1}$ instead
of $\I$ that
$\Th(R)$ is an $\Omega$-spectrum, and we let $\THH(R)$ be the $0$th
space.

In order to define the spectrum level multi-trace, we need to model
the additive structure of a spectrum in a very precise way. We next
explain how this can be done.

\subsection{The cyclic Barratt-Eccles construction}
Let $E_{\bullet}\Sigma_n$ be the cyclic set
$[k]\mapsto\Sigma_n^{k+1}$ with simplicial operators
\begin{align*}
&d_i(\sigma_0,\dots,\sigma_k)
=(\sigma_0,\dots,\sigma_{i-1},\sigma_{i+1},\dots,\sigma_k),\quad 0\leq
i\leq k,\\
&s_i(\sigma_0,\dots,\sigma_k)
=(\sigma_0,\dots,\sigma_{i-1},\sigma_i,\sigma_{i},\dots,\sigma_k),\quad
0\leq i\leq k,
\end{align*}
and cyclic operator
$t_k(\alpha_0,\dots,\alpha_k)
=(\alpha_k,\alpha_0,\dots,\alpha_{k-1})$. We let $\mathbb
E^{\infty}$ be the cyclic Barratt-Eccles operad with $n$th
space $E_{\bullet}\Sigma_n$, see \cite{BE},
\cite[6.5]{May}. This is an
$E_{\infty}$ operad in the sense that the realization $E\Sigma_n$
of the $n$th space is $\Sigma_n$-free and contractible.
We use the
notation $\mathbb E^{\infty}_{\bullet}$ for the associated functor
from based spaces to simplicial based spaces,
\begin{equation}\label{Einftyfunctor}
\mathbb E^{\infty}_{\bullet}(X)=\bigg(\coprod_{n\geq 0}
E_{\bullet}\Sigma_n\times X^n\bigg)/\sim,
\end{equation}
where the equivalence relation $\sim$ is defined as follows.
Notice first that 
the correspondence $\mathbf n\mapsto E_{\bullet}\Sigma_n$ defines a
contravariant functor from $\I$ to simplicial sets: Given a morphism
$\alpha\co\mathbf m\to \mathbf n$ in $\I$ and $\sigma\in \Sigma_n$, the
composition $\sigma\alpha$ has a unique factorization
$\sigma\alpha=\sigma_*(\alpha) \alpha^*(\sigma)$ with
$\sigma_*(\alpha)\co\mathbf m\to\mathbf n$ injective and order
preserving and $\alpha^*(\sigma)\in \Sigma_m$. In this way $\alpha$
induces a simplicial map,
$$
\alpha^*\co E_{\bullet}\Sigma_n\to E_{\bullet}\Sigma_m, \quad
(\sigma_0,\dots,\sigma_k)\mapsto
(\alpha^*(\sigma_0),\dots,\alpha^*(\sigma_k))
$$
and given $\beta\co\mathbf l\to \mathbf m$ it is clear that
$(\alpha\beta)^*=\beta^*\alpha^* $. Secondly, given a based space
$X$ the correspondence $\mathbf n\mapsto X^n$ defines a covariant
functor on $\I$ by letting a morphism $\alpha\co\mathbf m\to \mathbf
n$ act on $\mathbf x\in X^m$ by $\alpha_*(\mathbf x)=\mathbf y$,
where
$$
y_j=\begin{cases} x_i,&\text{ if }\alpha(i)=j\\
*,&\text{ if } j\notin \alpha(\mathbf m).
\end{cases}
$$
With this notation the equivalence relation in (\ref{Einftyfunctor})
is generated by the relations
$$
(e,\alpha_*(\mathbf x))\sim (\alpha^*(e),\mathbf x)
\quad\text{for } e\in E_k\Sigma_n,\ \mathbf x\in X^m\text{ and }
\alpha\co\mathbf m\to\mathbf n.
$$
In other words, $\mathbb E^{\infty}_{\bullet}(X)$ 
is the tensor product of the
functors $\mathbf n\mapsto E_{\bullet}\Sigma_n$ and $\mathbf n\mapsto X^n$
over $\I$, ie, the coend of the $\I^{\text{op}}\times \I$-diagram 
$E\Sigma_m\times X^n$, cf. \cite[IX.6]{ML}.  
We let $\mathbb E^{\infty}(X)$ be the
realization.
(Barratt and Eccles use the notation $\Gamma^+(X)$, but we want to
avoid this since we also use $\Gamma$-spaces in the sense of
Segal.)
We write the elements of $\mathbb E^{\infty}_{\bullet}(X)$ as
$[\mathbf \sigma,\mathbf x]$ where $\mathbf \sigma\in
E_{\bullet}\Sigma_k$ and $\mathbf x\in X^k$. Block sums of
permutations give $\mathbb E^{\infty}_{\bullet}(X)$ the structure of
a simplicial topological monoid,
$$
[\mathbf \sigma,\mathbf x]\cdot [\mathbf \sigma',\mathbf
  x']=[\mathbf\sigma\oplus\mathbf \sigma',(\mathbf x,\mathbf x')].
$$
The homotopy theoretical significance of the functor  $\mathbb
E^{\infty}(X)$ is that it provides a combinatorial model of 
$\Omega^{\infty}\Sigma^{\infty}(X)$ for non-degenerately based 
connected $X$. In more detail,
it is proved in \cite{BE} that in the diagram
$$
\mathbb E^{\infty}(X)\to\colim \Omega^n\mathbb
E^{\infty}(S^n\wedge X)\la\colim \Omega^n(S^n\wedge X),
$$
the left hand arrow is an equivalence for connected $X$ and the right
hand arrow is an equivalence in general.

We extend $\mathbb E^{\infty}$ to a functor on (symmetric) 
spectra by applying it in each spectrum dimension, ie, $\mathbb
E^{\infty}(E)_n=\mathbb E^{\infty}(E_n)$ with structure maps 
$$
S^1\wedge \mathbb E^{\infty}(E_n)\to \mathbb E^{\infty}
(S^1\wedge E_n)\to \mathbb E^{\infty}(E_{n+1}). 
$$ 
Since we assume spectra to be connective and convergent,
it easily follows that the natural map $E\to \mathbb E^{\infty}(E)$ is an
equivalence. Similarly, given a simplicial spectrum,  
we may apply $\mathbb E^{\infty}_{\bullet}$ degree-wise to get a bisimplicial
spectrum and then restrict to the simplicial diagonal. This is in
effect what we shall do when defining the spectrum level multi-trace.  

\subsection{The spectrum level multi-trace}
The multi-trace for a symmetric ring spectrum $R$ is a natural map
of multi-indexed spectra
\begin{equation}\label{spectrummultitrace}
\tr\co\underbrace{M_n(R)\wedge\dots\wedge M_n(R)}_{k+1}\to \mathbb
E^{\infty}_k(\underbrace{R\wedge\dots\wedge R}_{k+1}).
\end{equation}
Let us first explain how to define this when $R$ is a spectrum of
based simplicial sets. In this case $\tr$ is based on a natural
transformation
$$
\tr\co M_n(X_0)\wedge\dots\wedge M_n(X_k)\to \mathbb
E^{\infty}_k(X_0\wedge\dots\wedge X_k),
$$
where $X_0,\dots,X_k$ are based sets and $M_n(X_i)=\Map_*(\mathbf
n_+,\mathbf n_+\wedge X_i)$. Suppose given an element $(A^0,\dots,
A^k)$ in the domain and use matrix notation to write
$A^i=(x^i_{s,t})$. Let $D$ be the  set of multi-indices
corresponding to the non-trivial summands in the multi-trace
formula, ie,
$$
D=\{(s_0,\dots,s_k)\co x^0_{s_k,s_0}\neq *,\dots,
x^k_{s_{k-1},s_k}\neq *\}.
$$
Since by definition the matrices have at most one non-base point
entry in each column, the projections $(s_0,\dots,s_k)\mapsto s_i$
give rise to injective maps $p_i\co D\to \mathbf n$. Suppose that $D$
has cardinality $m$ and order the elements by choosing a bijection
$\gamma\co\mathbf m\to D$. The composition $p_i\gamma\co\mathbf m\to
\mathbf n$ is injective for each $i$ and admits a unique
factorization $p_i\gamma=\alpha_i\sigma_i$, where $\alpha_i$ is
injective and order
preserving and $\sigma_i\in \Sigma_m$. Consider the natural map
$$
D\to X_0\wedge\dots \wedge X_k,\quad (s_0,\dots,s_k)\mapsto
(x^0_{s_k,s_0},\dots,x^k_{s_{k-1},s_k})
$$
and let $\mathbf x$ be the composition
$$
\mathbf x\co\mathbf m\to D\to X_0\wedge\dots\wedge X_k.
$$
The first observation is that the element
$$
[(\sigma_0,\dots,\sigma_k);\mathbf x]\in
\Sigma_m^{k+1}\times_{\Sigma_m} (X_0\wedge\dots\wedge
X_k)^m
$$
is independent of the ordering $\gamma $ used to define it. By
definition the multi-trace is the image in
$\mathbb E^{\infty}_k(X_0\wedge\dots\wedge X_k)$,
\begin{equation}\label{multi-traceequation}
\tr(A^0,\dots,A^k)=[(\sigma_0,\dots,\sigma_k);\mathbf x]\in
\mathbb E^{\infty}_k(X_0\wedge\dots\wedge X_k).
\end{equation}
The second observation is that because of the base point relations
in the target this construction is natural with respect to based
maps in $X_0,\dots X_k$.
\begin{example}
As an example to illustrate the construction we calculate
$$
\tr\left(
\begin{bmatrix}
*&x^0_{12}\\ x^0_{21}&*
\end{bmatrix},
\begin{bmatrix}
*&x^1_{12}\\ x^1_{21}&*
\end{bmatrix}\right)
=[(1_2,\tau);(x^0_{21},x^1_{12}),(x^0_{12},x^1_{21})],
$$
where $\tau\in\Sigma_2$ is the non-identity element.
\end{example}

The spectrum level multi-trace (\ref{spectrummultitrace}) is defined
by degree-wise extending the above natural transformation to a natural
transformation between functors of simplicial sets.
We then extend this to a natural
transformation of  multi-indexed spectra by applying it in each
multi-degree. This gives the required maps
$$
\tr\co M_n(R_{n_0})\wedge\dots\wedge M_n(R_{n_k})\to \mathbb
E^{\infty}_k(R_{n_0}\wedge\dots\wedge R_{n_k}).
$$
In the case where $R$ is a spectrum of topological spaces we
observe that the expression in (\ref{multi-traceequation}) also
makes sense if $X_0,\dots,X_k$ are (non-degenerately) based
topological spaces, and we define $\tr$ by the same formula.

\subsection{The topological trace map}
We define a combinatorially enriched version $\Th^+(R)$ of
topological Hoch\-schild homology by applying B\"okstedt's
construction to the multi-indexed spectrum on the right hand side of
(\ref{spectrummultitrace}), ie,
$$
\Th_{\bullet}^+(R,n)\co[k]\mapsto
\Omega^{\infty}(\mathbb E^{\infty}_k(\underbrace{R\wedge\dots\wedge R}_{k+1}
\wedge S^n)).
$$
This is in a natural way the cyclic diagonal of a bicyclic spectrum. 
The $+$ decoration indicates that $\Th_{\bullet}^+(R,n)$ is a
homotopy commutative cyclic monoid. 
Using the natural inclusion $X\mapsto \mathbb
E_{\bullet}^{\infty}(X)$ we get a degree-wise equivalence
$\Th_{\bullet}(R)\to\Th_{\bullet}^+(R)$ and thus an equivalence of
realizations $\Th(R)\stackrel{\sim}{\to}\Th^+(R)$. The spectrum
level multi-trace has formal properties similar to the linear
multi-trace and in particular there results a cyclic map
$$
\tr\co\Th_{\bullet}(M_n(R))\to\Th_{\bullet}^+(R).
$$
(One can show that the realization can be extended to give an 
equivalence of $\mathbb T$-equivariant
spectra, but we shall not use this here.)
The definition of the topological trace map is now completely
analogous to the linear case. There is an obvious embedding of
cyclic spaces $B^{\cy}_{\bullet}GL_n(R)\to
\THH_{\bullet}(M_n(R))$ induced by the natural transformation
\begin{equation}
\begin{aligned}\label{BcyTHHinclusion}
&\Omega^{n_0}(M_n(R_{n_0}))^*\times \dots\times
\Omega^{n_k}(M_n(R_{n_k}))^*\\
&\to \Omega^{n_0+\dots+
n_k}(M_n(R_{n_0})\wedge\dots\wedge M_n(R_{n_k}))
\end{aligned}
\end{equation}
that sends a tuple of maps to their smash product. Composing with
the multi-trace we get a cyclic map
$$
\coprod_{n\geq 0}B_{\bullet}^{\cy}GL_n(R) \to \coprod_{n\geq
0}\THH_{\bullet}(M_n(R)) \to\THH_{\bullet}^+(R).
$$
This is a monoid homomorphism with respect to block sums of matrices
on the domain and the simplicial monoid structure on the target.
After realization and group completion we get maps
$$
\K^{\cy}(R)\to \Omega B(\THH^+(R))\stackrel{\sim}{\la}\THH^+(R)
\stackrel{\sim}{\la}\THH(R).
$$
The topological trace map
$\tr\co\K(R)\to\THH(R)$ is the homotopy class obtained by composing
with the homotopy section $\K(R)\to \K^{\cy}(R)$.

\begin{remark}
It is not difficult to extend this definition of the trace map
to a map of spectra or to refine it to a version of the cyclotomic
trace $\text{trc}\co\K(R)\to\operatorname{TC}(R)$, cf. \cite{BHM}.
However, this is not the purpose of
the present paper. A construction of the
trace map from a more categorical point of view has been given by
Dundas and McCarthy \cite{DM} and Dundas \cite{D}.
\end{remark}

Letting $n=1$ in (\ref{BcyTHHinclusion}) gives a map $B^{\cy}GL_1(R)\to
\THH(R)$.
The following proposition is immediate from the definitions
\begin{proposition}\label{BcyTHHproposition}
There is a strictly commutative diagram of spaces
\begin{equation}\label{Bcyinclusiondiagram}
\begin{CD}
B^{\cy}GL_1(R) @>>> \THH(R)\\
@VVV @VV\sim V\\
\K^{\cy}(R) @>\tr >> \Omega B\THH^+(R).
\end{CD}
\end{equation}
\vspace{-30pt}

\qed
\end{proposition}

\section{$\Gamma$-spaces and units of commutative ring
spectra}\label{gammasection}
 In this section we show that if $R$
is a commutative (and convergent) ring spectrum, then $GL_1(R)$ is
the $0$th space of an $\Omega$-spectrum. The same is true for the
group-like monoid $X_{h\I}$ associated to a commutative and
group-like $\I$-monoid $X$,  and we formulate the construction
in this generality.
\subsection{$\Gamma$-spaces}
We first recall Segal's notion of
$\Gamma$-spaces and the Anderson-Segal method for constructing the
associated homology theory. The paper by Bousfield and Friedlander
\cite{BF} is the basic
reference for this material. Let $\Gamma^o$ denote the category of
finite pointed sets and pointed maps. A $\Gamma$-space is a
functor $A\co\Gamma^o\to \mathcal T$ such that $A(*)=*$. We say that
a $\Gamma$-space is \emph{special} if given pointed sets $S$ and
$T$ the natural map $A(S\vee T)\to A(S)\times A(T)$  is an
equivalence. This implies in particular that $A(S^0)$ has the
structure of a homotopy associative and commutative H-space with
multiplication
$$
A(S^0)\times A(S^0)\simeq A(S^0\vee S^0)\to A(S^0).
$$
We say that $A$ is \emph{very special} if $A(S^0)$ is group-like,
ie, if the monoid of components is a group. A $\Gamma$-space $A$
extends to a functor on the category of pointed simplicial sets
in a two stage procedure. First $A$
is extended to the category of all pointed sets by forcing it to
commute with colimits. Given a simplicial set $X$ we then apply
$A$ degree-wise to get a simplicial space $[k]\mapsto A(X_k)$ with
realization $A(X)$. The main result is
that if $A$ is very special then the resulting functor is a homology
theory: Applying $A$ to a cofibration sequence of pointed
simplicial sets $X\to Y\to Y/X$ gives a homotopy fibration sequence
$$
A(X)\to A(Y)\to A(Y/X)
$$
in the sense that the inclusion of $A(X)$ in the homotopy fiber of
the second map is an equivalence. In particular, a very special
$\Gamma$-space gives rise to a symmetric $\Omega$-spectrum
$\{A(S^n)\co n\geq 0\}$, in which
the structure maps are the realizations of
the obvious (multi)-simplicial maps $S^1_{\bullet}\wedge
A(S^n_{\bullet})\to A(S^{n+1}_{\bullet})$.

\subsection{$\Gamma$-spaces associated to commutative $\I$-monoids}
In order to motivate the construction we recall
the definition of the $\Gamma$-space associated to a
commutative topological monoid $G$. Given a finite based set $S$, let
$\bar S$ be the subset obtained by excluding the base point.
Then $G(S)=G^{\bar S}$, and a based map $\alpha\co S\to T$ induces a map
$G(S)\to G(T)$ by multiplying the elements in $G$ indexed by
$\alpha^{-1}\{t\}$ for each $t\in\bar T$.

Implementing this idea for a commutative $\I$-monoid requires some
preparation. Given $S$ as above, let $\mathcal
P(\bar S)$ be the category of subsets and inclusions in $\bar S$.
A based map $\alpha\co S\to T$ induces a functor $\alpha^*\co\mathcal
P(\bar T)\to \mathcal P(\bar S)$ by letting
$\alpha^*(U)=\alpha^{-1}(U)$ for $U\subseteq \bar T$. 
We define a category $\D(S)$ of $\bar S$-indexed sum diagrams in $\I$ as
follows. An object is a functor $\theta\co\mathcal P(\bar S)\to\I$ 
that takes disjoint unions to coproducts of finite sets, ie, 
if $U\cap V=\emptyset$, then the diagram $\theta_U\to\theta_{U\cup
V}\la \theta_V$ represents $\theta_{U\cup V}$ as a coproduct of
finite sets. (The
category $\I$ itself does of course not have coproducts.) 
Notice in particular that $\theta_{\emptyset}=\mathbf 0$. 
A morphisms in $\mathcal D(S)$ is a natural
transformations of functors (not necessarily an isomorphism). 
This construction is clearly functorial in $\Gamma^{o}$: A based map
$\alpha\co S\to T$ induces a functor $\alpha_*\co\D(S)\to \D(T)$ by letting
$\alpha_* \theta=\theta\circ\alpha^*$.
Notice that an object $\theta$ in $\mathcal D(S)$ is determined by 
\begin{itemize}
\item
its 
values $\theta_s$ for $s\in \bar S$; 
\item
a choice of injective map $\theta_s\to \theta_U$
whenever $s\in U$ such that the induced map
$
\sqcup_{s\in U}\theta_s\to\theta_U
$
(with any ordering of the summands) is a bijection.
\end{itemize}
  
Restricting to the one-point
subsets of $\bar S$ gives a functor $\pi_S\co\D(S)\to \I^{\bar S}$,
where the latter denotes the product category indexed by $\bar S$ (we
let $\I^{\emptyset}$ denote the one-point category). This is an
equivalence of categories, and specifying an ordering of $\bar S$
gives a canonical choice of an inverse functor $\I^{\bar S}\to \D(S)$,
using the monoidal structure of $\I$. Notice however, that $\I^{\bar S}$ is
not functorial in $\Gamma^{o}$ as is the case for $\mathcal D(S)$.
\begin{lemma}\label{Dequivalencelemma}
Given a functor $Y\co\I^{\bar S}\to
\mathcal T$, the natural map
$$
\hocolim_{\D(S)}Y\circ \pi_S\to
\hocolim_{\I^{\bar S}}Y
$$
induced by $\pi_S$ is an equivalence.
\end{lemma}
\begin{proof}
By the cofinality criterion in \cite[XI.9.2]{BK}
(or rather its dual version) it
suffices to check that for any object $\mathbf a\in \I^{\bar S}$, the
category $(\mathbf a\downarrow \pi_S)$ of objects under $\mathbf a$
is contractible. But this is clear since this category has an
initial object.
\end{proof}

Let now $X$ be a commutative $\I$-monoid and consider the 
$\I^{\bar S}$-diagram $X^{\bar S}$ defined by
$$
\{\mathbf n_s\co s\in \bar S\}\mapsto \prod_{s\in\bar S}X(n_s).
$$
For $\bar S=\emptyset$ this should be interpreted as the
one-point space. We use the notation $X(S)$
for the $\D(S)$-diagram obtained by composing with $\pi_S$. With this
definition $X(S)$ is functorial in $S$ is the sense that a based
map $\alpha\co S\to T$ gives rise to a natural transformation of
$\mathcal D(S)$ diagrams
$$
X(S)\to X(T)\circ \alpha_*.
$$
In order to see this fix an object $\theta$ in $\mathcal D(S)$ and
choose an ordering of the subsets $U_t=\alpha^{-1}(t)$ for each $t\in
\bar T$. The
map in question is then a product over $\bar T$ of maps of the
form
$$
\prod_{s\in U_t}X(\theta_s)\to
X(\bigsqcup_{s\in U_t}\theta_s)\to 
X(\theta_{U_t}),
$$
where the first arrow comes from the multiplication in $X$ and the
second arrow is induced by the bijection $\sqcup_{s\in
U_t}\theta_s\to \theta_{U_t}$ determined by the sum diagram $\theta$.   
The main point is that since $X$ is commutative the composite map does
not depend on the ordering of $U_t$ used to define it.

By definition the $\Gamma$-space associated to $X$ is given by
$$
X_{h\I}(S)=\hocolim_{\mathcal D(S)}X(S).
$$
Given a based map $\alpha\co S\to T$, the induced map
$\alpha_*\co X_{h\I}(S)\to X_{h\I}(T)$ is the composition
$$
\hocolim_{\D(S)}X(S)\to \hocolim_{\D(T)}X(T)\circ
\alpha_*\to\hocolim_{\D(T)}X(T),
$$
where the first map is induced by the above natural transformation
and the second map is the map of homotopy colimits determined by the
functor $\alpha_*$.

It follows immediately from the definition that $X_{h\I}(S^0)=X_{h\I}$.
In order to compare $X_{h\I}(S^1)$ to the usual bar construction of
$X_{h\I}$ we specify an ordering of the $k$-simplices in $S^1_{\bullet}$
by letting
\begin{equation}\label{simplexordering}
u_j=(\underbrace{0,\dots,0}_j,1,\dots,1), \quad\text{for
}j=0,\dots,k.
\end{equation}
Then $u_0$ is the base point and $\bar S^1_k=\{u_1,\dots u_k\}$.
\begin{proposition}\label{BXproposition}
The $\Gamma$-space associated to a commutative $\I$-monoid $X$
is always special and is very special if and only if the
underlying monoid $X_{h\I}$ is group-like. In general
there is a natural equivalence $BX_{h\I}\stackrel{\sim}{\to}
X_{h\I}(S^1)$.
\end{proposition}
\begin{proof}
Using Lemma \ref{Dequivalencelemma} we get an equivalence
$$
X_{h\I}(S)=\hocolim_{\D(S)}X(S)\stackrel{\sim}{\to}\hocolim_{\I^{\bar
    S}}X^{\bar S} \cong
\prod_{s\in \bar S}X_{h\I},
$$
which is the condition for $X_{h\I}$ to be special. The statement about
being very special follows from the definition. In order
to define the equivalence we use the ordering of the
simplices of $S^1_{\bullet}$ given by (\ref{simplexordering}). As
noted earlier this ordering determines an equivalence
$$
\hocolim_{\I^k}X(n_1)\times \dots\times X(n_k)\to
\hocolim_{\mathcal D(S^1_k)}X(S^1_k)
$$
using the monoidal structure of $\I$. Identifying the left hand side
with the $k$-simplices of $B_{\bullet}X_{h\I}$ we get a simplicial map
$B_{\bullet}X_{h\I}\to X_{h\I}(S^1_{\bullet})$. Since this is an equivalence in
each simplicial degree its realization is also an equivalence as required.
\end{proof}

\begin{remark}
This construction of $\Gamma$-spaces based on the category $\mathcal
D(S)$ differs from that of Segal \cite[\S2]{Se} in that we
allow \emph{all} natural transformations, not only the
natural isomorphisms. Consequently, our definition of the
$\Gamma$-space associated to a commutative $\I$-monoid takes into account
all the maps $X(m)\to X(n)$ induced by morphisms in $\I$. As an
example, consider the commutative $\I$-monoid given by the classifying
spaces $BO(n)$ of the orthogonal groups. In this case Segal's
construction produces a special
$\Gamma$-space with underlying space $\coprod BO(n)$, whereas our
construction produces a very special $\Gamma$-space with
underlying space $\hocolim_{\I}BO(n)\simeq BO$. The last equivalence
follows from B\"okstedt's Lemma \ref{approximationlemma}.
Thus the two constructions respectively produce models of the
$(-1)$-connected and 0-connected topological K-theory spectrum.
\end{remark}

\begin{definition}\label{GammaGL1R}
Given a commutative (and convergent) symmetric ring spectrum $R$, let $GL_1(R)$
be the $\Gamma$-space associated to the $\I$-monoid $\Omega^n(R_n)^*$
considered 
in Section \ref{unitsection}, and let $gl_1(R)$ be the associated spectrum.
\end{definition}
It will always be clear from the context whether $GL_1(R)$ denotes a
$\Gamma$-space 
as above or the underlying group-like monoid as in Section \ref{unitssection}.

\section{Commutative ring spectra and splittings}\label{splittingsection}
Let $R$ be a commutative symmetric ring spectrum. In this section we
show that the 
natural inclusions $GL_1(R)\to B^{\cy}GL_1(R)$ and
$\Omega^{\infty}(R)\to \THH(R)$ have compatible left inverses in the
homotopy category, where by compatible we mean that these splittings are
related by a homotopy commutative diagram
\begin{equation}
\label{homotopysplitting}
\begin{CD}
B^{\cy}GL_1(R) @>r>> GL_1(R) \\
@VVV  @VVV\\
\THH(R) @>r>> \Omega^{\infty}(R).
\end{CD}
\end{equation}
We then define $\eta_R$ to be the composite homotopy class
$$
\eta_R\co BGL_1(R)\to K(R) \stackrel{\tr}{\to} 
\THH(R)\stackrel{r}{\to} \Omega^{\infty}(R).
$$
Using the diagrams (\ref{Bcyinclusiondiagram}) and
(\ref{homotopysplitting}) we get an alternative description as follows.
\begin{proposition}\label{firstpartproposition}
The homotopy class $\eta_R$ is represented by the composition
$$
BGL_1(R)\to L(BGL_1(R)) \stackrel{\sim}{\leftarrow}B^{\cy}GL_1(R)
\stackrel{r}{\to} GL_1(R),
$$
(where the first map is the inclusion of the constant loops), followed
by the inclusion of $GL_1(R)$ in $\Omega^{\infty}(R)$.
\qed
\end{proposition}
This concludes the first part of the proof of Theorem
\ref{introductiontheorem}. 
In order to motivate the construction, consider the cyclic bar construction
of a commutative monoid $G$. In this case the inclusion $G\to B^{\cy}(G)$ is
split by degree-wise multiplication in $G$. This can also be expressed in
terms of the $\Gamma$-space associated to $G$: The sequence
$G\to B^{\cy}G\to BG$ is the effect of evaluating $G$ on the
cofibration sequence 
$S^0\to S^1_{\bullet+}\to S^1_{\bullet}$, and the splitting is induced by the
projection $S^1_{\bullet+}\to S^0$ that maps $S^1_{\bullet}$ to the non-base
point in $S^0$. Let now $GL_1(R)$ be the
$\Gamma$-space 
defined in Definition \ref{GammaGL1R}. The next Lemma shows that we may replace
$B^{\cy}GL_1(R)$ by $GL_1(R)(S^1_+)$ up to homotopy.
\begin{proposition}\label{Bcyproposition}
Let $X$ be a convergent and commutative $\I$-monoid. Then there exists a space
$W^{\cy}$ and equivalences
$$
B^{\cy}X_{h\I}\stackrel{\sim}{\to} W^{\cy}\stackrel{\sim}{\la} X_{h\I}(S^1_+).
$$
\end{proposition}
\begin{proof}
 Let $W^{\cy}_{\bullet}$ be the cyclic space
$$
[k]\mapsto\hocolim_{\I^{k+1}\times \mathcal
  D(S^1_{k+})}X(n_0\sqcup\theta_{0})\times \dots\times X(n_k\sqcup
\theta_{k}),
$$
where $\theta$ is an object of $\mathcal D(S^1_{k+})$ and we write
$\theta_i=\theta_{u_i}$. We have functors $\I^{k+1}\to
\I^{k+1}\times\mathcal D(S^1_{k+})\la \mathcal D(S^1_{k+})$
obtained by fixing the initial object in one of the factors. The
cyclic structure of $W^{\cy}_{\bullet}$ is the obvious one such that the
induced maps
$$
B^{\cy}_{\bullet}X_{h\I}\to W^{\cy}_{\bullet}\la X_{h\I}(S^1_{\bullet+})
$$
become maps of cyclic spaces. Since we assume that $X$ is convergent,
it follows from B\"okstedt's approximation lemma
\ref{approximationlemma} (with $\I^{k+1}$ instead of $\I$), that these
maps are equivalences in each simplicial degree. After realization we
thus get a pair of equivalences relating $B^{\cy}X_{h\I}$ and
$X_{h\I}(S^1_+)$.
\end{proof}
\begin{remark}
The condition that $X$ be convergent is necessary for the argument
in Proposition \ref{Bcyproposition}, since otherwise the map
$$
\hocolim_{\I}X(m)\to \hocolim_{\I\times \I}X(m\sqcup n)
$$
induced by the functor $I\to I^2$, $\mathbf m\mapsto (\mathbf
m,\mathbf 0)$ need not be an equivalence. The $\I$-monoid $X(n)=X^n$
considered in Section \ref{tracesection} provides a counter example.
It should also be noted that the construction of the simplicial map
$B_{\bullet}X_{h\I} \to X_{hI}(S^1_{\bullet})$ in the proof of
Proposition \ref{BXproposition} cannot be applied to give a cyclic
map $B^{\cy}_{\bullet}X_{h\I}\to X_{h\I}(S^1_{\bullet+})$.
\end{remark}
Using the above equivalences, we define the splitting $r$ to be the
composite homotopy class
$$
r\co B^{\cy}GL_1(R)\simeq GL_1(R)(S^1_+)\to GL_1(R),
$$
where the last map is induced by the projection
$S^1_{\bullet+}\to S^0$.

We next consider a version $\Th'(R)$ of topological Hochschild
homology that relates to $\Th(R)$ as $GL_1(R)(S^1_+)$ relates to
$B^{\cy}GL_1(R)$. By definition this is the realization of the
cyclic spectrum
$$
\Th'_{\bullet}(R,n)\co[k]\mapsto \hocolim_{\mathcal D(S^1_{k+})}
\Omega^{\theta_0\sqcup\dots\sqcup\theta_k}(R_{\theta_0}\wedge\dots\wedge
R_{\theta_k}\wedge S^n),
$$
where $\theta$ denotes an object in $\mathcal D(S^1_{k+})$ and we again
write $\theta_i=\theta_{u_i}$. The cyclic structure maps are
defined as for $GL_1(R)(S^1_{\bullet+})$. For example, in spectrum
degree zero, $d_0\co\THH'_1(R)\to\THH'_0(R)$ is the composition
$$
\hocolim_{\theta\in\mathcal
D(S^1_{1+})}\Omega^{\theta_0\sqcup\theta_1}(R_{\theta_0}\wedge
R_{\theta_1})\to \hocolim_{\theta\in\mathcal
D(S^1_{1+})}\Omega^{\theta_{01}}(R_{\theta_{01}})\to
\hocolim_{\psi\in\mathcal D(S^1_{0+})}\Omega^{\psi_0}(R_{\psi_0}).
$$
Here $\theta_0\to\theta_{01}\la\theta_1$ denotes an object in
$\mathcal D(S^1_{1+})$ and the first map is induced by the natural
transformation that takes
$f\in\Omega^{\theta_0\sqcup\theta_1}(R_{\theta_0}\wedge
R_{\theta_1})$ to the element in
$\Omega^{\theta_{01}}(R_{\theta_{01}})$ given by the composition
$$
S^{\theta_{01}}\stackrel{\sigma^{-1}}{\lra}S^{\theta_0\sqcup\theta_1}=
S^{\theta_0}\wedge
S^{\theta_1}\stackrel{f}{\lra}R_{\theta_0}\wedge
R_{\theta_1}\stackrel{\mu_{\theta_0,\theta_1}}{\lra}
R_{\theta_0\sqcup\theta_1}\stackrel{\sigma}{\lra}R_{\theta_{01}},
$$
where $\sigma\co\theta_0\sqcup\theta_1\to\theta_{01}$ is the
bijection determined by $\theta$. The second map is induced by the
natural transformation $d_0\co\mathcal D(S^1_{1+})\to \mathcal
D(S^1_{0+})$. With this definition we have the equality
$$
d_0=d_1\co\Th'_1(R)\to \Th'_0(R)
$$
and consequently the iterated boundary maps give a well-defined cyclic map
$r\co\Th'_{\bullet}(R)\to \Th'_0(R)$, where the target is considered
a constant cyclic spectrum. In spectrum degree zero we thus get a
cyclic map of spaces
$$
r\co\THH'_{\bullet}(R)\to\Omega^{\infty}(R).
$$
The next proposition is the analogue of Proposition
\ref{Bcyproposition}.
\begin{proposition}\label{TH'proposition}
The spectra $\Th(R)$ and $\Th'(R)$ are related by a pair of equivalences.
\end{proposition}
\begin{proof}
Letting $W^{\cy}_{\bullet}=\{W^{\cy}_{\bullet}(n)\co n\geq 0\}$ 
denote the cyclic spectrum
$$
[k]\mapsto\hocolim_{\I^{k+1}\times \mathcal D(S^1_{k+})}
\Omega^{\mathbf n_0\sqcup\theta_0\sqcup\dots\sqcup \mathbf
n_k\sqcup\theta_k} (R_{\mathbf n_0\sqcup\theta_0}\wedge\dots\wedge
R_{\mathbf n_k\sqcup\theta_k}\wedge S^n),
$$
the proof proceeds exactly like the proof of Proposition
\ref{Bcyproposition}.
\end{proof}

As in the definition of the trace map we consider the
transformation
$$
\Omega^{\theta_0}(R_{\theta_0})^*\times\dots\times
\Omega(R_{\theta_k})^*\to
\Omega^{\theta_0\sqcup\dots\sqcup\theta_k}(R_{\theta_0}\wedge\dots\wedge
R_{\theta_k})
$$
that sends a tuple of maps to their smash product. Viewing these
maps as natural transformations of $\mathcal D(S^1_{k+})$-diagrams
we get a cyclic map
$$
GL_1(R)(S^1_{\bullet+})\to\THH'_{\bullet}(R).
$$
It follows immediately from the definitions that this map is
compatible with the splittings of $B^{\cy}GL_1(R)$ and $\THH'(R)$
in the sense of the following proposition.
\begin{proposition}
There is a strictly commutative diagram of spaces
$$
\begin{CD}
GL_1(R)(S^1_+)@>>>GL_1(R)\\
@VVV @VVV \\
\THH'(R)@>>> \Omega^{\infty}(R).
\end{CD}
$$
\vspace{-30pt}

\qed
\end{proposition}

The homotopy commutative diagram (\ref{homotopysplitting}) in the
beginning of this section is derived from this using the
equivalences in Proposition \ref{Bcyproposition} and Proposition
\ref{TH'proposition}.

\section{The Hopf map and free loops on infinite loop
  spaces}\label{Hopfsection}
In this section we finish the proof of Theorem \ref{introductiontheorem}
by showing that the composite homotopy class
$$
BGL_1(R)\to
L(BGL_1(R))\stackrel{\sim}{\leftarrow}B^{\cy}GL_1(R)\stackrel{r}{\to}GL_1(R) 
$$
is multiplication by $\eta$ in the sense explained in the introduction.
More generally, let $G$ be a very special $\Gamma$-space and let
$g=\{G(S^n)\co n\geq 0\}$ be the associated $\Omega$-spectrum.
By \cite[4.1]{BF},  the $\Gamma$-space $G$ is determined by $g$ in the sense
that the diagram
$$
G(X)\stackrel{\sim}{\to}\hocolim\Omega^nG(S^n\wedge X) 
\stackrel{\sim}{\leftarrow}\hocolim\Omega^n(G(S^n)\wedge
X))=\Omega^{\infty}(g\wedge X)
$$
specifies a natural equivalence $G(X)\simeq\Omega^{\infty}(g\wedge
X)$.
Evaluating $G$ on the based cyclic
set $S^1_{\bullet+}$ gives a cyclic space $G(S^1_{\bullet+})$. The
realization $G(S^1_+)$ then has a $\mathbb T$-action and, as in the case
of the cyclic bar construction, we consider the composite map
$$
\mathbb T\times G(S^1_+)\to G(S^1_+)\to G(S^1)
$$
with adjoint $G(S^1_+)\to L(G(S^1))$. In the next proposition we
analyze the homotopy fibration sequence obtained by evaluating $G$ on
the cofibration sequence $S^0\to S^1_{\bullet+}\to S^1_{\bullet}$.
\begin{proposition}\label{gammafreeloopproposition}
There is a commutative diagram of homotopy fibration sequences
\begin{equation}\label{gammafreeloopdiagram}
\begin{CD}
G(S^0) @>>> G(S^1_+) @>>> G(S^1)\\
@VV\sim V @VV\sim V @ |\\
\Omega(G(S^1))@>>> L(G(S^1))@>>> G(S^1)
\end{CD}
\end{equation}
in which the the vertical maps are equivalences.
\end{proposition}
\begin{proof}
The commutativity of the right hand square is immediate since we
evaluate a loop at the unit element of $\mathbb T$. In order to
prove commutativity of the left hand square we recall that for any
cyclic space $X_{\bullet}$, the $\mathbb T$-action on the zero
simplices $X_0\subseteq |X_{\bullet}|$ has the following
description. If $x$ is an element of $X_0$ and $u\in \mathbb T$,
$$
u\cdot x=[t_1s_0x,u]\in |X_{\bullet}|.
$$
Here $t_1$ is the cyclic operator in degree one and we make the
identification $\mathbb T=\Delta^1/\partial\Delta^1$. Using this, it
is easy to check that the composition
$$
\mathbb T\times G(S^0)\to \mathbb T\times G(S^1_+)\to G(S^1_+)\to
G(S^1)
$$
is given by $(u,x)\mapsto [x',u]$, where $x'\in G(S^1_1)$ is the
image of $x$ under the homeomorphism induced by the based
bijection $S^0\to S^1_1$. The above composition clearly equals the
composition
$$
\mathbb T\times G(S^0)\to S^1\wedge G(S^0)\to G(S^1),
$$
which shows that the left hand square in the diagram is also
commutative.
\end{proof}

In Diagram (\ref{gammafreeloopdiagram}) the map $L(G(S^1_+))\to G(S^1)$
in the lower sequence is split by the inclusion of the constant loops,
and the map $G(S^0)\to G(S^1_+)$ in the upper sequence is split by
evaluating $G$ on the projection $r\co S^1_{\bullet+}\to S^0$ that maps
$S^1_{\bullet}$ to the non-base point of $S^0$. The next proposition
expresses the fact that these splittings are not compatible in general.
As usual $\eta\in \pi_1^s(S^0)$ denotes the stable Hopf map.
\begin{proposition}\label{etasplittingproposition}
Using the natural equivalences $\Omega^{\infty}(g\wedge S^1)\simeq
G(S^1)$ and $\Omega^{\infty}(g)\simeq G(S^0)$, the composite homotopy class
$$
G(S^1)\to L(G(S^1))\stackrel{\sim}{\la}G(S^1_+)\stackrel{r}{\to}G(S^0),
$$
is given by $\Omega^{\infty}(g\wedge \eta)$.
\end{proposition}
\begin{proof}
We first observe that (\ref{gammafreeloopdiagram}) is in
fact a diagram of infinite loop spaces and infinite loop maps, and
that as such it is equivalent to the following diagram of spectra
$$
\begin{CD}
g\wedge S^0 @>>> g\wedge S^1_+ @>>>g\wedge
S^1\\
@V\simeq VV @V\simeq V\Phi V @|\\
 F(S^1,g\wedge S^1)@>>> F(S^1_+,g\wedge S^1)@>>> F(S^0,g\wedge S^1).
\end{CD}
$$
Here $F(-,g\wedge S^1)$ is the obvious function spectrum and
the upper and lower cofibration sequences are both induced from
$S^0\to S^1_+\to S^1$. These cofibration sequences have canonical
stable splittings induced by the projection $r\co S^1_+\to S^0$ and
the associated stable section $s\co S^1\to S^1_+$. 
The vertical map $\Phi$
in the middle is the adjoint of
$$
S^1_+\wedge g\wedge S^1_+ \to g\wedge S^1_+ \to g\wedge S^1,
$$
where the first map uses the action of $S^1$ on itself given by the
group structure. It is clear that the above diagram is equivalent to
the one obtained by smashing $g$ with
$$
\begin{CD}
S^0 @>>> S^1_+ @>>> S^1\\
@V\simeq VV @V\simeq V\Phi V @| \\
F(S^1,S^1)@>>> F(S^1_+,S^1)@>>> F(S^0,S^1).
\end{CD}
$$
Here the definition of $\Phi$ is analogous to the
definition given above. We must prove that the stable map
$$
S^1=F(S^0,S^1)\stackrel{r^*}{\to} 
F(S^1_+,S^1)\stackrel{\Phi}{\simeq} S^1_+\stackrel{r}{\to}S^0
$$
represents $\eta$. Using the canonical splittings to represent
$\Phi^{-1}$ as a $2\times 2$ matrix, this composition represents the
off diagonal term. It is therefore the negative of the composite
stable map starting in the upper right corner of the diagram,
$$
S^1\stackrel{s}{\to} S^1_+ \stackrel{\Phi}{\to} F(S^1_+,S^1)
\stackrel{s^*}{\to}F(S^1,S^1).
$$
The adjoint of this is the stable map
$$
S^1\wedge S^1 \stackrel{s\wedge s}{\to} 
S^1_+\wedge S^1_+ \to S^1_+\to S^1,
$$
where the second map is the group multiplication in $S^1$. It
is well-known that this composition represents $\eta$.
For example, one can see this by considering the equivariant splitting
$\widetilde \Sigma S^1\to S^1\wedge S^1_+$, whose domain is the
unreduced suspension 
of $S^1$, and then use that the map of homotopy colimits induced by
the diagram 
$$
\left(*\la S^1\times S^1\to *\right )\to \left(*\la S^1\to *\right)
$$
represents a generator of $\pi_3(S^2)$, cf. \cite[XI.4]{Wh}. The
result now follows since $\eta$ has order two.
\end{proof}
\noindent\textbf{Proof of Theorem \ref{introductiontheorem}}\qua
The only thing left to prove is that the homotopy class
$BGL_1(R)\to GL_1(R)$ considered in Proposition
\ref{firstpartproposition} agrees with the one
in Proposition
\ref{etasplittingproposition} when the $\Gamma$-space in question is
$GL_1(R)$ and we use the canonical equivalence
$BGL_1(R)\stackrel{\sim}{\to} GL_1(R)(S^1)$, cf. Proposition \ref{BXproposition}.
Let $W$ be the realization of the simplicial space
$$
[k]\mapsto\hocolim_{\I^k\times \mathcal D(S^1_k)}X(\mathbf
n_1\sqcup \theta_1)\times\dots\times X(\mathbf n_k\sqcup
\theta_k).
$$
Since the space $W^{\cy}$ considered in Proposition \ref{Bcyproposition}
is the realization of a cyclic space it has a
$\mathbb T$-action, and the adjoint of the composition
$$
\mathbb T\times W^{\cy}\to W^{\cy}\to W
$$
gives an equivalence $W^{\cy}\to L(W)$. It easily follows that we have a
commutative diagram of equivalences
$$
\begin{CD}
B^{cy}GL_1(R)@>\sim>>W^{cy}@<\sim <<GL_1(R)(S^1_+)\\
@VV \sim V @ VV\sim V  @VV \sim V \\
L(BGL_1(R))@>\sim>> L(W)@<\sim<< L(GL_1(R)(S^1)).
\end{CD}
$$
It thus suffices to check that the homotopy class defined by the diagram
$$
BGL_1(R)\stackrel{\sim}{\to}W\stackrel{\sim}{\la}GL_1(R)(S^1)
$$
is compatible with the canonical equivalence
$BGL_1(R)\stackrel{\sim}{\to}GL_1(R)(S^1)$. We do
this by exhibiting an explicit homotopy inverse of the equivalence
$GL_1(R)(S^1)$ $\to W$. Let $W'$ be the realization of the simplicial
space
$$
[k]\mapsto \hocolim_{\mathcal D(S^1_k)\times\mathcal D(S^1_k)}
X(\psi_1\sqcup\theta_1)\times \dots\times X(\psi_k\sqcup
\theta_k).
$$
Ordering the $k$-simplices $S^1_k$ as in the proof of
\ref{BXproposition} gives an equivalence
$W\stackrel{\sim}{\to}W'$. Moreover, using the monoidal structure
of $\I$ we get a functor $\mathcal D(S^1_k)\times \mathcal
D(S^1_k)\to \mathcal D(S^1_k)$. Varying $k$ this is a
transformation of simplicial categories and consequently there is
an induced simplicial map $W_{\bullet}'\to GL_1(R)(S^1_{\bullet})$. It
is easy to  check that the composition $W\to W'\to GL_1(R)(S^1)$ is a
homotopy inverse of the map in question and that the composition with
$BGL_1(R)\to W$ is the canonical equivalence. This completes the proof. \qed

\begin{remark}
The discussion in the above proof can be generalized to any
commutative $\I$-monoid $X$. If
$X$ is convergent, $B^{\cy}X_{h\I}$ and $X_{h\I}(S^1_+)$ are related by
cyclic equivalences as in Proposition \ref{Bcyproposition}, and if
$X$ is group-like we get equivalences relating the same spaces by
comparing them to the relevant free loop spaces. In case $X$ is
both convergent and group-like, the two equivalences agree by an
argument similar to the above.
\end{remark}

{\bf Acknowledgement}\qua{The author was partially supported by a
grant from the NSF.}

\end{document}